\documentclass[3p,authoryear]{elsarticle}

\usepackage{amsmath,amssymb,amsthm}
\usepackage[colorlinks=true,linkcolor=blue,citecolor=blue,urlcolor=blue]{hyperref}

\theoremstyle{plain}
\newtheorem{theorem}{Theorem}

\newtheorem{proposition}{Proposition}
\newtheorem{corollary}{Corollary}

\theoremstyle{remark}

\newcommand{\R}{\mathbb{R}}
\newcommand{\PR}{\mathbf{P}}
\newcommand{\E}{\mathbf{E}}
\DeclareMathOperator{\var}{var}

\newcommand{\boldi}{{\bf i}}
\newcommand{\Inm}{{I_{n,m}}}

\newcommand{\secref}[1]{Section~\ref{sec:#1}}

\newcommand{\propref}[1]{Proposition~\ref{prop:#1}}

\newcommand{\thmref}[1]{Theorem~\ref{thm:#1}}

\newcommand{\corollaryref}[1]{Corollary~\ref{cor:#1}}

\journal{Journal of Statistical Planning and Inference}

\ifpreprint\fi 

\begin{document}

\begin{frontmatter}

	\title{Sharp  upper and lower bounds on the variance of U-statistics\tnoteref{t1}}

	\begin{abstract}
		For a  U-statistic with a  symmetric kernel  of any  rank, we offer  new upper and lower  bounds 
		on its variance, which are uniformly sharper than those established in the classical paper by  \citet[Theorem 5.2]{Hoeffding1948}. Our new bounds have strong optimality properties and can  tightly describe how the variance scales with the degree of the kernel.   As an application, we   leverage our new upper bound to improve an existing Berry-Esseen theorem for the normal approximation of incomplete U-statistics based on Bernoulli sampling, in a computationally frugal regime where  the computational budget is much smaller than the available sample size.
		

	\end{abstract}
	
	\author[1]{Kiarash Nikoo}
\ead{knikoo@student.unimelb.edu.au}
	
\author[2]{Liqian Zhang}
\ead{liqian.zhang@postgrad.manchester.ac.uk}
	
\author[1]{Dennis Leung\corref{cor1}} 
\ead{dennis.leung@unimelb.edu.au}
\cortext[cor1]{Corresponding author}

	\affiliation[1]{organization={School of Mathematics and Statistics,
University of Melbourne},
city={Parkville},
state={VIC 3031},
country={Australia}
	}
\affiliation[2]{organization={Department of Mathematics,   University of Manchester},
addressline={Oxford Road},
city={Manchester},
postcode={M13 9PL},
country={the United Kingdom}
}
 

\end{frontmatter}

\section{Introduction}\label{sec:intro}

 Let $X_{1},\dots,X_{n}$ be independent and identically distributed sample points taking values in a   measurable space $(M,\mathcal M)$, and  $h\colon M^{m}\to\R$ be a measurable real-valued function  symmetric in its $m$ arguments, i.e. $h(x_1, \dots, x_m) = h(x_{\pi_1}, \dots, x_{\pi_m})$  for any permutation $(\pi_1, \dots, \pi_m)$ of the indices $(1, \dots, m)$. For  $m \leq n$,
in  his seminal paper, \citet{Hoeffding1948} introduced the U-statistic
\begin{equation} \label{U_stat_def}
U_{n}\equiv\binom{n}{m}^{-1}\sum_{ 1 \leq i_1 < \cdots <  i_m \leq n}h(X_{i_1}, \cdots, X_{i_m}) 
\end{equation}
that averages    the evaluations of the kernel function $h$ over all  possible ${n \choose m}$ subsets of $m$ distinct sample points. Throughout, for any natural number $k \leq n$, we will adopt the shorthand $[k] \equiv \{1, \dots, k\}$; moreover,  given  a function $t:M^k \rightarrow \mathbb{R}$ in $k$  arguments in $M$ and another natural number $\ell$, we will use $\E[t]$, $\var[t]$,  and $\E[t^\ell]$  as respective shorthands for $\E[t(X_1, \dots, X_k)]$, $\var[t(X_1, \dots, X_k)]$ and $\E[(t(X_1, \dots, X_k))^\ell]$.  We will assume that $0 < \var[h] < \infty$,  so that $U_n$ is not  always trivially equal to its expectation $\theta \equiv \E[h]$.

 

The  structure of the variance of $U_n$ is  well understood. For any $c \in [m]$,   let 
\begin{equation} \label{Psi_c_def}
\Psi_c (x_1, \dots, x_c) \equiv \E[ h(x_1, \dots, x_c, X_{c+1}, \dots, X_m)]- \theta  
\end{equation}
be the centered conditional expectation of $h$ by fixing $c$ of its $m$ arguments, and 
\begin{equation} \label{g_c_def}
g_c(x_1, \dots, x_c) \equiv \Psi_c (x_1, \dots, x_c)  - \sum_{j = 1}^{c-1} \sum_{1 \leq i_1 < \dots < i_j\leq c} g_j (x_{i_1}, \cdots, x_{i_j})  
\end{equation}
be the recursively defined \emph{canonical projection function} of $h$.
Define the second moment quantities 
\begin{equation}\label{second_moment_quantities}
\zeta_{c}\equiv\E [ \Psi_c^2   ] \quad  \text{ and } \quad \delta_c \equiv \E[g_c^2]  
\end{equation}
for every $c \in [m]$. 
 The two exact variance formulas
\begin{align}
\var[U_{n}] &= {n \choose m}^{-1} \sum_{c = 1}^m {m \choose c} {n - m \choose m - c} \zeta_c  \label{exact_variance_formula_zeta}\\
  &= \sum_{c = 1}^m {m \choose c}^2 {n \choose c}^{-1} \delta_c  \label{exact_variance_formula_delta}
\end{align}
are  standard results in the U-statistic literature  \citep[Chapter 1]{KB1994Ustat}. 
We also remark that $\zeta_m=\var[h]$ by definition.

Nevertheless,  instead of the exact forms in \eqref{exact_variance_formula_zeta} and \eqref{exact_variance_formula_delta},  it is sometimes more instructive to work with \emph{bounds} on $\var[U_{n}]$, particularly when one demands  a more transparent understanding of how $\var[U_{n}]$ scales with the kernel degree $m$. In this article, we will contribute new bounds on $\var[U_{n}]$ that are sharp, which we believe will be useful for   future research on the theory of U-statistics,   an area that has proven  to be  fertile even after nearly 8 decades since the   work by  \citet{Hoeffding1948}.


\section{Main results}
 \citet[Theorem 5.2]{Hoeffding1948} established that, if the kernel $h$ is of \emph{rank} at least $d$, then 
\begin{equation}\label{hoeff}
	\frac{m}{d} K_{n, m}(d)\zeta_{d} \leq \var[U_{n}]  \leq  K_{n, m}(d) \zeta_m 
\end{equation}
for the constant
\begin{equation} \label{K_const_def}
 K_{n, m}(d)  \equiv  \binom{n}{m}^{-1}\sum_{c=d}^{m}\binom{m-1}{c-1}\binom{n-m}{m-c}.
\end{equation}
The rank of a kernel $h$ (or of the corresponding U-statistic)  is defined to be the smallest integer $r \in [m]$ such that   
$
\zeta_1  = \cdots = \zeta_{r-1} = 0$ and $\zeta_r >0$; it  is well-defined, because  it must be true that $0 \leq \zeta_1 \leq \cdots \leq \zeta_m$  as  an easy consequence of Jensen's inequality.  Alternatively, for any $c \in [m]$,  using the well-known identities   that
\begin{equation}\label{zeta_c_in_deltas}
\zeta_c =  \sum_{j = 1}^c \binom{ c}{  j}\delta_j 
\end{equation}
and 
\begin{equation}\label{delta_c_in_zetas}
\delta_c = \sum_{j = 1}^c (-1)^{c - j} \binom{c}{  j} \zeta_j
\end{equation}
the rank can be equivalently  defined as the smallest  $r \in [m]$ such that   
$
\delta_1  = \cdots = \delta_{r-1} = 0$  and $\delta_r >0$.\footnote{Standard references, such as \citet[p.32]{KB1994Ustat}, may equivalently write  \eqref{zeta_c_in_deltas} and \eqref{delta_c_in_zetas} as $\zeta_c =  \sum_{j = 1}^c \binom{ c}{c - j}\delta_j $ and $\delta_c = \sum_{j = 1}^c (-1)^{c - j} \binom{c}{c - j} \zeta_j$ using the symmetry of binomial coefficients.} In the literature, a kernel of rank $1$ is said to be \emph{non-degenerate}, and a kernel of rank $m$ is said to be \emph{completely degenerate}.

Rather surprisingly, in this article we    illustrate that the classical inequalities in  \eqref{hoeff} are  actually sub-optimal. Our main contributions  are the following   upper and lower bounds   on $\var[U_{n}]$ that are  \emph{uniformly sharper} than \eqref{hoeff}, which are proven in \ref{sec:pf_main}.

\begin{theorem}[Sharper variance bounds on U-statistics]\label{thm:sharp}
	Suppose $U_n$ is a U-statistic defined as  in \eqref{U_stat_def}   of rank at least $d \in [m]$ and $m < n$.
	Then
	\begin{equation}\label{sharp}
		\binom{m}{d} B_{n,m}(d) \zeta_{d}  \leq \var[U_{n}]  \leq B_{n,m}(d) \zeta_m,
\end{equation}
	for the constant
	\begin{equation} \label{B_mnd_def}
		B_{n,m}(d) \equiv \binom{n}{d}^{-1}\binom{m}{d} .
	\end{equation}
	Moreover, these two properties hold:
	\begin{enumerate}[(i)]
	\item It is true that 
			\begin{equation}\label{KB_intra_bound_nice_looking}
		\frac{m}{d} K_{n, m}(d) \leq \binom{m}{d} B_{n,m}(d)	 \quad \text{ and }	\quad  B_{n,m}(d) \leq   K_{n, m}(d)  ,
			\end{equation}
			and  each inequality in \eqref{KB_intra_bound_nice_looking} becomes an equality if and only if $d\in\{1,m\}$.
		\item Each inequality in \eqref{sharp} is an equality if and only if $\delta_{d+1} = \cdots = \delta_m = 0$.

	\end{enumerate}

\end{theorem}

We assume $m <n$ in \thmref{sharp} to exclude the mathematically uninteresting case $m = n$, in which the variance upper bounds in both \eqref{hoeff}  and \eqref{sharp}   trivially become equalities as $B_{m,m}(d)= K_{m,m}(d) = 1$ and $\var[U_n]= \zeta_m$. Property $(i)$ in \thmref{sharp} says that   our   bounds in  \eqref{sharp} are strictly sharper than those in  \eqref{hoeff}    by having better constants, except when the kernel is non-degenerate  or completely degenerate; it is established through a neat argument that expresses both $B_{n,m}(d)$ and $K_{n,m}(d)$ in terms of the same hypergeometric random variable and exploits its binomial moments; see \eqref{Kmnd_in_hypergeom} and \eqref{Bmnd_in_hypergeom}. 

To further illustrate the optimality of our result, property $(ii)$ in \thmref{sharp}   says that both inequalities in \eqref{sharp} can be \emph{attained} when   $\delta_d$ is the only non-vanishing quantity among all of $\delta_1, \dots, \delta_m$, that is,    when $U_n$ is    entirely comprised of the $d$-th component in its canonical  decomposition \citep[p.24]{KB1994Ustat} and therefore, upon centering by its mean $\theta$, has the form
\begin{equation} \label{only_dth_component}
U_n - \theta = {m \choose d} {n \choose d}^{-1} \sum_{1 \leq  {i_1} < \cdots <  {i_d}\leq n} g_d(X_{i_1} ,  \cdots,  X_{i_d}).
\end{equation}
An example of a kernel that can give rise to \eqref{only_dth_component}  is as follows. Let $\psi\colon M\to\R$ be a (measurable) function such that $\E[\psi]=0$ and $0<\E[\psi^{2}]<\infty$, and define the mean-zero kernel
	\begin{equation}\label{extremal}
		h_{\psi}(x_{1},\dots,x_{m})\equiv\sum_{1\le i_{1}<\cdots<i_{d}\le m}
		\psi(x_{i_{1}})\cdots\psi(x_{i_{d}}).
	\end{equation}
	For instance, if $A\in\mathcal M$ is a set such that $0<\PR(X_{1}\in A)<1$  and   $\psi(x)\equiv\mathbf{1}\{x\in A\}-\PR(X_{1}\in A)$, then   $\E[\psi]=0$ and $\E[\psi^{2}]=\PR(X_{1}\in A)\,\PR(X_{1}\notin A)\in(0,\infty)$. 
In \ref{sec:show_only_d_component}, we will show that a U-statistic based on the kernel $h_{\psi}$ above satisfies \eqref{only_dth_component}.

To the best of our knowledge, \thmref{sharp} is   new. That being said,  a prototype of our upper bound on $\var[U_{n}]$ in \eqref{sharp} has been hinted by   \citet{Peng2022}, who  set out to improve upon the Berry-Esseen theorem  for a non-degenerate $U_n$ stated in \citet[Theorem 10.3]{ChenGoldsteinShao2011}. In this latter line of work, a key step is to bound the variance of the \emph{H\'{a}jek projection remainder}  of $U_n$
\[
U_{\bar{h}_m} \equiv {n \choose m}^{-1} \sum_{1 \leq i_1 < \cdots < i_m \leq n} \bar{h}_m (X_{i_1}, \dots, X_{i_m}),
\]
 which is   itself a U-statistic   constructed with the kernel $\bar{h}_m (x_1, \dots, x_m) \equiv  h (x_1, \dots, x_m) - \theta - \sum_{i = 1}^m \Psi_1 (x_i)$ and is   of rank at least $2$. \citet[ App.~C.2, p.~279]{Peng2022}    applied a bound   qualitatively comparable to our upper bound   in \eqref{sharp} to get  
 \begin{equation} \label{peng_bound}
 \var[U_{\bar{h}_m}] \leq \frac{m^2}{n^2} \E[\bar{h}_m^2].\footnote{Under $m < n$, our upper bound in \eqref{sharp} at $d=2$ actually gives the slightly better bound   $\var[U_{\bar{h}_m}] \leq \frac{m(m-1)}{n(n-1)}\E[\bar{h}_m^2 ]$.}
 \end{equation}
 On the contrary,   \citet{ChenGoldsteinShao2011}  have effectively  applied  \citet{Hoeffding1948}'s upper bound in \eqref{hoeff}  at $d=2$ on $\var[U_{\bar{h}_m}]$, and their subsequent arguments \citep[p.284-286]{ChenGoldsteinShao2011}, which show the inequality
  \begin{equation} \label{chen_shao_bdd_on_K}
  K_{n, m}(2) \leq \frac{m(m-1)^2}{n (n - m+1)},
  \end{equation}
   will further lead to the bound $\var[U_{\bar{h}_m}] \leq \frac{m (m-1)^2}{n (n - m+1)} \E[\bar{h}_m^2]$, which is fundamentally looser than \eqref{peng_bound} with respect to order in $m$.  As a result, among other things,  the Berry-Esseen bound for $U_n$ in \citet[Theorem 3]{Peng2022} is superior to the one in \citet[Theorem 10.3]{ChenGoldsteinShao2011} by exhibiting a tighter dependence on the degree $m$. The proof of our upper bound in \eqref{sharp} essentially  sharpens \citet{Peng2022}'s  arguments   and generalizes them for a U-statistic of any arbitrary rank. In \secref{BE_ICU}, we will further leverage this new upper bound to improve an existing Berry-Esseen theorem for  incomplete U-statistics with Bernoulli sampling.
 
 We have also developed the proposition below (proven in \ref{sec:proof_proposition})  to illustrate some properties of  the constant $B_{n, m}(d)$, particularly how it scales with the kernel degree $m$:

\begin{proposition}[Further properties of the constant $B_{n, m}(d)$]\label{prop:K}
  The following inequalities are true:
	\begin{enumerate}
		\item[(i)] For any $d \in [m]$,
		 
			\begin{equation}\label{Bpowers}
				\left(\frac{m - d +1}{ n}\right)^{d}\le B_{n, m}(d) \le\left(\frac{m}{n}\right)^{d}.
			\end{equation}
		
		\item[(ii)] Suppose $n > 2m$. For any $1 < d < m$,
	
			\begin{equation}\label{Klower}
			 \frac{ (n - m)_{m - d}}{ (n - d)_{m - d}}	 \binom{m-1}{d-1}B_{n, m}(d)    \leq K_{n,m}(d),
			\end{equation}
			where $(a)_b$ denotes the falling factorial $a (a-1)\cdots (a- b +1) $.
	\end{enumerate}
\end{proposition}

\propref{K}$(i)$ should intuitively make sense to the reader once we rewrite the constant as 
\begin{equation} \label{rewrite_B_as_product}
B_{n, m}(d) = \prod_{j=0}^{d-1}\frac{m-j}{n-j};
\end{equation}
it suggests that for a fixed $d$, $B_{n, m}(d) $ is of the order $O(m^d /n^d)$. To understand \propref{K}$(ii)$, first note that it is always true that $K_{n,m}(d) \leq \binom{m-1}{d-1}B_{n, m}(d)$, which simply rephrases the first inequality in \eqref{KB_intra_bound_nice_looking}. Therefore, \eqref{Klower} can be understood as a ``reverse'' of the latter inequality via adjustment by the factor $\frac{ (n - m)_{m - d}}{ (n - d)_{m - d}} < 1$. However, when $d$ is held fixed and $m^2/n \longrightarrow 0$, the factor  $ \frac{ (n - m)_{m - d}}{ (n - d)_{m - d}} $ is very close to  $1$, which suggests that $K_{n,m}(d)$ is approximately larger than $B_{n,m} (d)$ by a factor of $ \binom{m-1}{d-1} = O(m^{d-1})$. Indeed, it is no coincidence that 
only an upper  bound of the order $O(m^3/n^2)$ can be established for $K_{n,m}(2)$ as in \eqref{chen_shao_bdd_on_K} by  \citet{ChenGoldsteinShao2011}, which is looser  than the tighter bound of the order $O(m^2/n^2)$ suggested by \eqref{Bpowers} at $d =2$.

\section{Application to   Berry-Esseen theorems for incomplete U-statistics with a small Bernoulli sampling budget}\label{sec:BE_ICU}
Let $I_{n, m} \equiv \{(i_1, \dots, i_m): 1 \leq i_1 < \cdots < i_m \leq n \}$, and for a positive integer $N < {n \choose m}$, define $p \equiv  N {n \choose m}^{-1}$. Suppose 
 $(Z_{\mathbf{i}})_{\mathbf{i}\in I_{n,m}}$ are ${n \choose m}$ independent and identically distributed Bernoulli random variables such that $\PR(Z_\mathbf{i} = 1)  = p$, and  that they are also independent of all the sample points $X_1, \dots, X_n$. Let $\hat{N} \equiv \sum_{\boldi \in \Inm} Z_\boldi$ be the count of how many of these Bernoulli random variables realize as a ``1''. For a kernel with a large degree $m$,  instead of $U_n$ in \eqref{U_stat_def} which sums over all ${n \choose m}$ kernel evaluations, to achieve computational savings one can construct the \emph{incomplete $U$-statistic with Bernoulli sampling} 
 \begin{equation} \label{icu_def}
	U'_{n,N}\equiv\frac{1}{\hat N}\sum_{\mathbf{i}\in I_{n,m}}Z_{\mathbf{i}}\,
	h(X_{\mathbf{i}}) ,
\end{equation}
where   we have used the shorthand $h(X_{\mathbf{i}}) \equiv h(X_{i_1}, \dots, X_{i_m}) $ for any $\boldi \equiv (i_1, \dots, i_m) \in \Inm$; by convention, $U'_{n,N}$ is set to be zero when $\hat{N} =0$. Because   only approximately a fraction $p$ of    all the Bernoulli samplers $Z_{\mathbf{i}}$ in \eqref{icu_def} will be ``switched on'', $U'_{n,N}$   sums over many fewer kernel evaluations than $U_n$, particularly when the \emph{computational budget} $N$ is small.

The normal convergence behavior of $U'_{n,N}$   has been well-studied. In particular, if $r$ is the rank of the kernel $h$, it has been established that, for a fixed $m$,
\begin{equation} \label{n_ll_N_convergence}
\frac{\sqrt{N}(U_{n, N}' - \theta)}{\zeta_m^{1/2}} \text{ converges in distribution to the standard normal as }n, N \longrightarrow \infty \text{ with } n^r/N \longrightarrow \infty;
\end{equation} 
see  \citet[Theorem 1.1 $(ii)$]{Leung2026} which has its roots in \citet[Corollary 1]{janson1984asymptotic}. Note that the asymptotic regime here  suggests a   computationally frugal situation where $N$ could be much less than $n^r$.  \citet[Theorem 3.2]{Leung2026}  also developed a Berry-Esseen theorem to gauge the normal approximating accuracy in terms of Kolmogorov distance  for the weak convergence  in \eqref{n_ll_N_convergence}. Turns out,  the said Berry-Esseen theorem  can be further improved with our new upper bound on $\var[U_{n}]$  in \eqref{sharp}. Below, we will first state the improved theorem, and then explain how it has improved upon its predecessor \citet[Theorem 3.2]{Leung2026}; we will use $\Phi(\cdot)$  to denote the standard normal  distribution function.

%
%
%

\begin{corollary}[Improved Berry-Esseen theorem for the normal convergence in \eqref{n_ll_N_convergence}] \label{cor:BE}
Let  $U_{n ,N}'$ be constructed as in \eqref{icu_def} with a symmetric kernel $h$ of rank $r$ on $M^m$. Suppose that  
\begin{equation} \label{h_mean_zero_assumption}
\E[h] = 0,
\end{equation}
$0 <  \E[|h|^4] < \infty$  and  $2 \leq m < n/2$ hold.  Then, there exists an absolute constant $C >0$ such that
\begin{multline}  \label{BE_bdd_Nlln_4th_moment}
 \sup_{z \in \mathbb{R}}\bigg| \PR\bigg(\frac{\sqrt{N} U_{n, N}'}{\zeta_m^{1/2}} \leq z\bigg) - \Phi(z) \bigg|
 \leq \\
 C \bigg\{ 
 \frac{\E[|h|^3] (1 - 2p + 2p^2)}{ \sqrt{N(1-p)}\zeta_m^{3/2}}  + \exp \bigg( \frac{- n \zeta_m^3}{24 m(\E[|h|^3])^2} \bigg)  + 
 \frac{  \sqrt{NB_{n,m}(r)}}{\sqrt{1-p}} +
\frac{\sqrt{m \E[(h^2 - \zeta_m)^2]}}{ \sqrt{n}\zeta_m}
\bigg\}.
\end{multline}
\end{corollary}

As is customary in the literature on Berry-Esseen theorems for U-statistics, we have made the mean zero assumption in \eqref{h_mean_zero_assumption} only to avoid cluttering the notation in  our bound \eqref{BE_bdd_Nlln_4th_moment}; obviously, one could present the theorem in terms of the centered kernel $h(\cdot) - \theta$ without making this assumption. Moreover, similarly to \citet[Theorem 3.2]{Leung2026}, a version of the result with the milder finite third moment assumption $\E[|h|^3] < \infty$ can also be stated, but we have opted to only present the notationally cleaner result \eqref{BE_bdd_Nlln_4th_moment} under the finite fourth moment assumption  $\E[|h|^4] < \infty$ for the purpose of illustration. 

The original Berry-Esseen bound in \citet[Theorem 3.2, Eq. $(3.7)$]{Leung2026} is essentially the same as   \eqref{BE_bdd_Nlln_4th_moment} except that the quantity $K_{ n, m}(r)$ was in the place of $B_{n, m}(r)$;  indeed, as explained in \ref{sec:proof_improved_BE}, the proof of \corollaryref{BE}  is only a minor modification of that of \citet[Theorem 3.2]{Leung2026} by making use of the new upper bound in \eqref{sharp}.  As such, when $1 < r <m$, our Berry-Esseen bound \eqref{BE_bdd_Nlln_4th_moment} is now tighter, not only because $B_{n, m}(r)$ is smaller than  $K_{n, m}(r)$  as per property $(i)$ in \thmref{sharp}, but also because it is smaller than it by approximately a  factor of  order  $O((m-1)^{1-r })$ for a fixed $r$ as $m^2/n \longrightarrow 0$, as per \eqref{Klower}. The latter is an important attribute for modern applications of incomplete U-statistics such as uncertainty inference in random forests, where the degree $m$ corresponds to the subsample size for growing each tree, which is typically quite large; we refer the readers to \citet[Section 1.2]{Leung2026} which contains more related discussion.


\section*{Acknowledgments}

This research is partially funded by the Manchester-Melbourne-Toronto Research Fund. 
Dennis Leung would also like to thank Dr. Yuta Koike for some helpful discussion during the early stage of this project. 

 \section*{Declaration of generative AI and AI-assisted technologies in the manuscript preparation process.}

During the preparation of this work, the first author (K.N.) used Claude (Anthropic) to assist in refining and verifying certain mathematical arguments and calculations, while the last author (D.L.) used Codex (OpenAI) to assist with some final edits. The authors independently verified all mathematical statements and proofs,   and take full responsibility for the content of the published article.

\appendix 
\section{Proof of \thmref{sharp}} \label{sec:pf_main}

From \eqref{exact_variance_formula_delta}, we can  first write  the variance of $U_n$ as
\begin{align}
	\var[U_{n}] &=\sum_{j=d}^{m} B_{n,m}(j)  \binom{m}{j}\delta_j \label{var_decomp_1}\\
	& = \sum_{j=d+1}^{m} B_{n,m}(j)  \binom{m}{j}\delta_j + B_{n, m}(d)  \binom{m}{d}\zeta_d,
	\label{var_decomp_2}
	\end{align}
 where \eqref{var_decomp_1} is true because
  \begin{equation}\label{delta_below_d_are_zero}
 \delta_1 = \dots = \delta_{d-1} = 0, 
 \end{equation} 
 as the rank of $h$ is at least $d$, and \eqref{var_decomp_2} is true because
 \begin{equation}\label{zeta_d_equal_delta_d}
 \zeta_{d}=\delta_{d},
 \end{equation}
  which can be seen by taking $c = d$ in \eqref{zeta_c_in_deltas} and using \eqref{delta_below_d_are_zero}.

Under $m < n$,  $B_{n,m}(r) /B_{n,m}(r+1)= \frac{n - r}{m - r} > 1$ for any $r \in [m-1]$, which implies
 \begin{equation} \label{decreasing_B_constants}
 B_{n,m}(r) \text{ is strictly decreasing in } r \in [m]. 
 \end{equation}
 In light of \eqref{decreasing_B_constants}, one can deduce from 
 \eqref{var_decomp_1} that
 \begin{equation} \label{main_result_RHS}
 	\var[U_{n}] =\sum_{j=d}^{m} B_{n,m}(j)  \binom{m}{j}\delta_j  \leq  B_{n,m}(d) \sum_{j = d}^m \binom{m}{  j} \delta_j    = B_{n,m}(d) \zeta_m,
 \end{equation}
 where   the last equality in \eqref{main_result_RHS}  comes from \eqref{zeta_c_in_deltas},  \eqref{delta_below_d_are_zero} and \eqref{zeta_d_equal_delta_d}.
On the other hand, from \eqref{var_decomp_2}, one can immediately see that 
  \begin{equation} \label{main_result_LHS}
  \binom{m}{d} B_{n,m}(d) \zeta_{d}  \leq  \sum_{j=d+1}^{m} B_{n,m}(j)  \binom{m}{j}\delta_j + B_{n, m}(d)  \binom{m}{d}\zeta_d =  \var[U_{n}].
  \end{equation}
  Not only do the inequalities  \eqref{main_result_RHS} and \eqref{main_result_LHS} together imply \eqref{sharp}, but it is also easy to see that each of them precisely becomes an equality if and only if 
  $\delta_{d+1} = \cdots = \delta_m = 0$, which   proves the property $(ii)$.
 
%

It remains to show property $(i)$. 
			To that end, we will  exploit expressions of $B_{n, m}(d)$ and $K_{n, m}(d)$ in terms of hypergeometric random variables, as well as an elementary result on the binomial moments of hypergeometric distributions \citep[  Eq. (93)]{pitman2001random}: If $H_{T, S, t}$ is a $\mathrm{Hypergeometric}(T,S,t)$  random variable, i.e. it is the number of successes in $t$ draws without replacement, out of a finite population of size $T$ containing $S$ successes, then 
 \begin{equation} \label{binomial_moment_hypergeom}
 \E\Bigg[\binom{H_{T, S, t}}{r} \mathbf{1}\{H_{T, S, t}\ge r\} \Bigg] = \frac{\binom{S}{r} \binom{t}{r} }{  \binom{T}{r}}  \text{ for } r = 1,\dots, S. \footnote{This result is also posed as an exercise in the undergraduate text by  \citet[p.257, Ex. 37]{pitman1993probability}. } 
 \end{equation}

 We will first prove 
 \begin{equation} \label{Bnmd_leq_Kmnd}
 B_{n, m}(d)\le K_{n, m}(d),
 \end{equation}
as well as the fact that
 \begin{equation} \label{Bnmd_leq_Kmnd_when_exact}
 \eqref{Bnmd_leq_Kmnd} \text{ is an exact inequality if and only if } d \in \{1, m\}.
 \end{equation}
 	 From the definition of $K_{n,m}(d)$ in \eqref{K_const_def} and $\binom{m-1}{c-1}=\frac{c}{m}\binom{m}{c}$, we can first write
	\begin{equation} \label{Kmnd_in_hypergeom}
		K_{n, m}(d) 
		=  \sum_{c=d}^{m}\frac{c}{m} \frac{\binom{m}{c}\binom{n-m}{m-c}}{\binom{n}{m}}  =  \sum_{c=d}^{m}\frac{c}{m}\PR(H_{n, m, m}=c)\ = \frac{1}{m}\E \bigg[  H_{n, m, m} \mathbf{1}\{H_{n, m, m}\ge d\} \bigg].
	\end{equation}
	On the other hand, using the definition of $B_{n, m}(d)$ in \eqref{B_mnd_def} and  formula \eqref{binomial_moment_hypergeom}, we have
	\begin{equation}  \label{Bmnd_in_hypergeom}
	B_{n, m}(d) = \binom{m}{d}^{-1}\E\bigg[\binom{H_{n, m, m}}{d}   \mathbf{1}\{H_{n, m , m}\ge d\}  \bigg]  = \frac{1}{m}\E\bigg[  \frac{\binom{H_{n, m, m}- 1}{d-1}}{\binom{m-1}{d-1}}   H_{n, m, m}   \mathbf{1}\{H_{n, m , m}\ge d\}  \bigg]. \footnote{For convenience, we set {}$\binom{H_{n, m, m}- 1}{d-1} = 0$ when $H_{n, m, m} = 0$.}
	\end{equation}
	Since $H_{n, m, m}\le m$ by definition, it is always true that $\frac{\binom{H_{n, m, m}- 1}{d-1}}{\binom{m-1}{d-1}}  \leq 1$;  therefore,  it is always true that
	\begin{equation}\label{something}
	 \frac{\binom{H_{n, m, m}- 1}{d-1}}{\binom{m-1}{d-1}}  H_{n, m, m}    \mathbf{1}\{H_{n, m , m}\ge d\} \quad \leq \quad  H_{n, m, m}   \mathbf{1}\{H_{n, m , m}\ge d\} ,
	\end{equation}
	which proves \eqref{Bnmd_leq_Kmnd} by comparing the rightmost expectations in   \eqref{Kmnd_in_hypergeom} and \eqref{Bmnd_in_hypergeom}. Moreover,  it is easy to see that \eqref{something} becomes an equality when $d \in \{1, m\}$, so the ``if'' part of \eqref{Bnmd_leq_Kmnd_when_exact} is proved. When $1 < d < m$, on the event $\{H_{n, m, m} = m -1\}$ which has  a strictly positive probability, \eqref{something} becomes a strict inequality because  $\frac{\binom{H_{n, m, m}- 1}{d-1}}{\binom{m-1}{d-1}}  = \frac{\binom{m-2}{d-1}}{\binom{m-1}{d-1}} < 1$, which proves the ``only if'' part of \eqref{Bnmd_leq_Kmnd_when_exact} by comparing  \eqref{Kmnd_in_hypergeom} and \eqref{Bmnd_in_hypergeom} again.

	Next, we will prove $\frac{m}{d} K_{n, m}(d) \leq \binom{m}{d} B_{n,m}(d)	$, which is equivalent to
		\begin{equation} \label{Kmnd_leq_binom_const_Bnmd}
		K_{n, m}(d) \leq \binom{m-1}{d-1} B_{n, m}(d),
	\end{equation}
        as well as  the fact that 
        \begin{equation} \label{Kmnd_equal_binom_const_Bnmd}
       \text{ \eqref{Kmnd_leq_binom_const_Bnmd} is an exact equality if and only if } d \in \{1, m\}
           \end{equation}
           in an analogous manner. The  rightmost expression of $B_{n, m}(d)$ in \eqref{Bmnd_in_hypergeom}  readily  gives
	\begin{equation} \label{upper_Bmnd_in_hypergeom}
		\binom{m-1}{d-1} B_{n, m}(d) = \frac{1}{m}\E\bigg[\binom{H_{n, m, m}-1}{d-1} H_{n, m, m}\mathbf{1}\{H_{n, m, m}\ge d\}\bigg].
	\end{equation}
	 Since it is always true that $\binom{H_{n, m, m}-1}{d-1}\ge1$ on the event $\{H_{n, m, m}\ge d\}$, we have	\begin{equation} \label{upper_hypergeom_comparison}
		H_{n, m, m}\mathbf{1}\{H_{n, m, m}\ge d\} \quad \leq \quad \binom{H_{n, m, m}-1}{d-1} H_{n, m, m}\mathbf{1}\{H_{n, m, m}\ge d\},
	\end{equation}
	which proves \eqref{Kmnd_leq_binom_const_Bnmd} by comparing the  expectations in \eqref{Kmnd_in_hypergeom} and \eqref{upper_Bmnd_in_hypergeom}. Moreover, it is easy to see that \eqref{upper_hypergeom_comparison} becomes an equality when $d\in\{1,m\}$, so the ``if'' part of \eqref{Kmnd_equal_binom_const_Bnmd} is proved. When $1<d<m$, on the event $\{H_{n, m, m}=m\}$ which has a strictly positive probability, we have $\binom{H_{n, m, m}-1}{d-1}=\binom{m-1}{d-1}>1$, so \eqref{upper_hypergeom_comparison} becomes a strict inequality, which proves the ``only if'' part of \eqref{Kmnd_equal_binom_const_Bnmd} by comparing \eqref{Kmnd_in_hypergeom} and \eqref{upper_Bmnd_in_hypergeom} again.

Collecting \eqref{Bnmd_leq_Kmnd}, \eqref{Bnmd_leq_Kmnd_when_exact}, \eqref{Kmnd_leq_binom_const_Bnmd} and \eqref{Kmnd_equal_binom_const_Bnmd}  proves the property $(i)$, and we are done.
 


\section{Proof that the kernel $h_\psi$ in \eqref{extremal} satisfies \eqref{only_dth_component}} \label{sec:show_only_d_component}

  In this section, all the special functions in \eqref{Psi_c_def} and \eqref{g_c_def} and  moment quantities in \eqref{second_moment_quantities} mentioned  are defined with respect to the specific kernel $h_{\psi}$. Note that  $\delta_1 = \dots =  \delta_{d-1} = 0$  because 
  \begin{equation} \label{first_g_zero}
	g_1= \dots = g_{d-1} = 0
  \end{equation} as per their recursive definitions in \eqref{g_c_def}, because it is easy to see that $\Psi_c (x_1, \dots, x_c) = 0$ for the kernel $h_{\psi}$ and  any $c \leq d-1$. To show that  $\delta_{d+1} = \dots =  \delta_m = 0$, because of \eqref{first_g_zero}, one first notes that $\Psi_d(x_1,\dots,x_d) = g_d  (x_1,\dots,x_d) =  \psi (x_1 ) \cdots \psi (x_d )$, and so $\delta _d =   (\E[\psi^2])^d $.  Further, by  exploiting independence, via direct calculations it must be true that
	 \[
	\zeta_m =  \E[h_\psi^2] = \sum_{1\le i_{1}<\cdots<i_{d}\le m}
		\E[\psi^2(X_{i_{1}})] \cdots \E[\psi^2(X_{i_{d}})] = \binom{m}{d} (\E[\psi^2])^d = \binom{m}{d} \delta _d.
	 \]
	In light of the generic formula for $\zeta_m$ in \eqref{zeta_c_in_deltas}, the above display  is only true when all  the $\delta_{d+1}, \cdots, \delta_m$ vanish.

\section{Proof of Proposition~\ref{prop:K} } \label{sec:proof_proposition}

We first prove $(i)$. 
The upper bound in \eqref{Bpowers}  follows from \eqref{rewrite_B_as_product} since $\frac{m-j}{n-j} \leq \frac{m}{n}$ when $m \leq n$. Moreover, since any $j$ is at most $d - 1$ in \eqref{rewrite_B_as_product}, we have $\frac{m - d+1}{n} \leq \frac{m-j}{n-j}$, and the lower bound in \eqref{Bpowers}  follows. Now we prove $(ii)$. It suffices to show
  \begin{equation} \label{suffice_for_comparison_between_K_and_B}
  \frac{ (n - m)_{m - d}}{ (n - d)_{m - d}}\binom{m-1}{d-1}B_{n, m}(d) =  \binom{m-1}{d-1}\binom{n-m}{m-d}\Big/\binom{n}{m},
  \end{equation}
  because the term on the right hand side of \eqref{suffice_for_comparison_between_K_and_B}
	  is the term at $c=d$ of the sum defining $K_{n,m}(d)$ in \eqref{K_const_def}. Since $\binom{n}{m}\binom{m}{d}=\binom{n}{d}\binom{n-d}{m-d}$, we have 
	\begin{align*}
	\binom{m-1}{d-1}\binom{n-m}{m-d}\big/\binom{n}{m} 
	&=  \binom{m-1}{d-1}B_{n, m}(d) \frac{ {n \choose d}{n - m \choose m - d}}{ \binom{n}{m} {m \choose d}}\\ 
	 &= \binom{m-1}{d-1}B_{n, m}(d) \frac{  {n - m \choose m - d}}{   {n - d \choose m - d}} \\
	& = \binom{m-1}{d-1}B_{n, m}(d) \frac{ (n - m) \cdots ((n - m)- (m - d)+1)}{ (n - d) \cdots ((n - d)- (m - d)+1)},
	\end{align*}
and our proof of \eqref{suffice_for_comparison_between_K_and_B} is done.

%
%
%
 
 \section{Proof of \corollaryref{BE}} \label{sec:proof_improved_BE}
 
 The proof is  almost word for word the same as the proof of \citet[Theorem 3.2]{Leung2026}; one   simply has to  replace all the occurrences of $K_{n, m}(r)$ in the proof of \citet[Theorem 3.2]{Leung2026} with $B_{n, m}(r)$. For the interested reader who will explore the fine details in \citet{Leung2026}, we  briefly explain why such simple replacements suffice below.

 The proof of \citet[Theorem 3.2]{Leung2026} relies on a technical result stated in \citet[Lemma 5.5 $(i)$]{Leung2026}, whose proof in \citet[Section D.1.3, p. 34]{Leung2026} in turn uses  that $\E[U_n^2] \leq K_{n,m}(r) \zeta_m$, which is simply  \citet{Hoeffding1948}'s upper bound in \eqref{hoeff} under \eqref{h_mean_zero_assumption}. Obviously, by our new upper bound in \eqref{sharp}, it is also true that $\E[U_n^2] \leq B_{n,m}(r) \zeta_m$, so $K_{n,m}(r)$ can be supplanted. Another place where $K_{n, m}(r)$ shows up in the original proof of \citet[Theorem 3.2]{Leung2026} is at \citet[p. 19]{Leung2026}, where the inequality ${n \choose m}^{-1} \leq K_{n,m}(r)$ is employed. However, it is certainly also true that ${n \choose m}^{-1} = \frac{m \cdots  1 }{n \cdots (n- m+1)} \leq B_{n,m}(r) =   \frac{m \cdots (m- r+1) }{n \cdots (n- r+1)}$ , so $K_{n, m}(r)$ can be supplanted there as well.

\bibliographystyle{elsarticle-harv}
\bibliography{references}

\end{document}